\documentclass[11pt]{article}
\usepackage[utf8]{inputenc}

\usepackage[letterpaper,margin=1in]{geometry}

\usepackage{graphicx}

\usepackage{xcolor}
\definecolor{darkblue}{rgb}{0,0,0.5}
\definecolor{darkred}{rgb}{.8,0,0}
\definecolor{darkgreen}{rgb}{0,.5,0}
\usepackage[colorlinks=true,linkcolor=darkblue,citecolor=darkblue,urlcolor=darkblue]{hyperref}

\usepackage{tikz}
\usetikzlibrary{shapes, arrows, calc, positioning}
\tikzstyle{vertex}=[draw,thick,fill=white,circle,inner sep=2pt]

\usepackage{subcaption}  

\usepackage{amsthm,amsmath,amssymb,amsfonts} 

\theoremstyle{plain}
\newtheorem{theorem}{Theorem}
\newtheorem*{theorem*}{Theorem}
\newtheorem{lemma}[theorem]{Lemma}
\newtheorem*{lemma*}{Lemma}

\newtheorem*{proposition*}{Proposition}
\newtheorem{corollary}[theorem]{Corollary}
\newtheorem*{corollary*}{Corollary}

\newtheorem*{claim*}{Claim}

\newtheorem*{conjecture*}{Conjecture}
\newtheorem{question}[theorem]{Question}

\theoremstyle{definition}
\newtheorem{definition}[theorem]{Definition}
\newtheorem*{definition*}{Definition}

\theoremstyle{remark}

\newtheorem*{example*}{Example}

\newcommand{\alontarsi}[1]{\ensuremath{\mathrm{AT}(#1)}}
\newcommand{\corres}[1]{\ensuremath{\chi_{DP}(#1)}}

\newcommand{\mad}[1]{\ensuremath{\mathrm{mad}\!\left(#1\right)}}

\title{Relation between the correspondence chromatic number and the Alon--Tarsi number}
\author{Eric Culver\\
Department of Mathematical and Statistical Sciences\\
University of Colorado Denver\\
\texttt{eric.culver@ucdenver.edu}
\and
Stephen G. Hartke\thanks{Supported in part by a Collaboration Grant from the Simons Foundation (\#316262 to Stephen G. Hartke).}\\
Department of Mathematical and Statistical Sciences\\
University of Colorado Denver\\
\texttt{stephen.hartke@ucdenver.edu}
}
\date{June 10, 2022}

\begin{document}

\maketitle

\begin{abstract}
  We study the relation between the correspondence chromatic number and the Alon--Tarsi number, both upper bounds on the list chromatic number of a graph.
  There are many graphs with Alon--Tarsi number greater than the correspondence chromatic number.
  We present here a family of graphs with arbitrary Alon--Tarsi number, with correspondence chromatic number one larger.
  \medskip
  
  
\end{abstract}

\section{Introduction}

List coloring was introduced independently by Vizing~\cite{Vizing1976} and by Erd\H{o}s, Rubin, and Taylor~\cite{ErdosRubinTaylor1979} as a generalization of graph coloring.
The list chromatic number $\chi_\ell(G)$ of a graph $G$ (also known as the \emph{choosability} of $G$) is the least size of lists assigned to each vertex such that a proper vertex coloring can be chosen for any such lists.

Alon and Tarsi~\cite{AlonTarsi1992} used the polynomial method known as the Combinatorial Nullstellensatz to prove an upper bound on the list chromatic number.
For a graph $G$, the Alon--Tarsi number $\alontarsi{G}$ is the best possible upper bound that can be proved using this method.
The study of the Alon--Tarsi number as a graph parameter was initiated by Jensen and Toft~\cite{JensenToft1994}.

Dvo\v{r}\'{a}k and Postle~\cite{DvorakPostle2018} introduced correspondence coloring as a generalization of list coloring, where on each edge a different correspondence between colors can be given that describes which colors conflict.
The correspondence chromatic number $\corres{G}$ of a graph $G$ is thus an upper bound on the list chromatic number $\chi_\ell(G)$.

Both the correspondence chromatic number and the Alon--Tarsi number are upper bounds on the list chromatic number of a graph.  What is the relationship between these two parameters?
There are many examples of graphs for which the Alon--Tarsi number is larger than the correspondence chromatic number.
Following from a result of Bernshteyn, all regular triangle-free graphs with large enough degree satisfy this, and it fact the ratio can become arbitrarily large in this direction.

However, there are also examples of graphs where the correspondence chromatic number is one more than the Alon--Tarsi number, such as even cycles (see Figure~\ref{fig:C_4} for the example of a cycle on four vertices) and a planar bipartite graph constructed by Bernshteyn and Kostochka~\cite{Bernshteyn2019}.
These examples have Alon--Tarsi numbers equal to two and three, respectively, so it is possible the correspondence chromatic number can exceed the Alon--Tarsi number only when the Alon--Tarsi number is small.
We show that this is not the case by constructing an infinite family of graphs with arbitrarily high Alon--Tarsi number which have correspondence chromatic number one larger.
In contrast with the situation when the correspondence chromatic number is less than the Alon--Tarsi number, we show the ratio cannot become arbitrarily large by presenting an upper bound on the correspondence chromatic number in terms of the Alon--Tarsi number.

\section{Preliminaries}

In this section we present formal definitions for proper coloring, list coloring, correspondence coloring, and the Alon--Tarsi number.
All graphs considered in this paper are simple and finite.  See West~\cite{West2001} for graph theory terminology and notation not defined here.

For a graph $G$, the \emph{maximum average degree $\mad{G}$} is defined as the maximum average degree over all subgraphs of $G$, and can be computed by
\[
  \mad{G} = \max_{\text{subgraph $H$ of $G$}} \frac{\sum_{v\in V(H)} \deg_H(v)}{|V(H)|}
          = \max_{\text{subgraph $H$ of $G$}} \frac{2|E(H)|}{|V(H)|}.
\]
The \emph{degeneracy} of a graph $G$ is the least integer $k$ such that every subgraph $H$ of $G$ has a vertex whose degree in $H$ is at most $k$, and can be computed by
\[
  \text{degeneracy of $G$} = \max_{\text{subgraph $H$ of $G$}} \delta(H),
\]
where $\delta(H)$ is the minimum degree of the subgraph $H$.
Note that the maximum average degree is an upper bound on the degeneracy of a graph.
By considering greedy coloring, the degeneracy of $G$ plus one is an upper bound on all of the coloring parameters of $G$ considered in this paper: the chromatic number, the list chromatic number, the correspondence chromatic number, and the Alon--Tarsi number.

\paragraph{Proper coloring.}
For a positive integer $k$, a \emph{proper $k$-coloring} of a graph $G$ is a function that assigns one of $k$ colors to each of the vertices of $G$ so that adjacent vertices have different colors.
For conciseness, we will assume in this paper that all colorings are proper (meaning that adjacent vertices have different colors), and so we will henceforth drop the adjective ``proper''.
We say that $G$ is \emph{$k$-colorable} if it has a $k$-coloring.
The \emph{chromatic number $\chi(G)$} of a graph $G$ is the least $k$ such that $G$ is $k$-colorable.

\paragraph{List coloring.}
A \emph{list assignment} for a graph $G$ is a function $L$ that assigns to each vertex a list of colors.
An \emph{$L$-coloring} of $G$ is a proper coloring of $G$ such that for each vertex, the color assigned to it is chosen from its list.
A graph $G$ is \emph{$k$-choosable} if there exists an $L$-coloring of $G$ for every assignment $L$ of lists of size at least $k$ to the vertices of $G$.
The \emph{list chromatic number $\chi_\ell(G)$} of a graph $G$ is the least $k$ such that $G$ is $k$-choosable.

Every $k$-choosable graph is $k$-colorable.
However, the converse is known to not be true in general.
For example, the complete bipartite graphs $K_{n,n}$ have chromatic number $2$, but are not $2$-choosable when $n\ge 3$.

\paragraph{Correspondence coloring.}
A \emph{correspondence assignment} for $G$ consists of a list assignment $L$ and a function $C$ that to every edge $uv \in E(G)$ assigns a partial matching $C_{uv}$ between $\{u\} \times L(u)$ and $\{v\} \times L(v)$. The Cartesian product is used to distinguish the vertices of $C_{uv}$ in case the same color appears in both $L(u)$ and $L(v)$.
An \emph{$(L,C)$-coloring} of $G$ is a function $\phi$ that to each vertex $v \in V(G)$ assigns a color $\phi(v) \in L(v)$, such that for every $uv \in E(G)$, the vertices $(u,\phi(u))$ and $(v, \phi(v))$ are non-adjacent in $C_{uv}$. We say that $G$ is \emph{$(L,C)$-colorable} if such an $(L,C)$-coloring exists.

Correspondence coloring generalizes list coloring since if we set each correspondence $C_{uv}$ to match exactly the common colors of $L(u)$ and $L(v)$ for each $uv \in E(G)$, then an $(L,C)$-coloring is an $L$-coloring.

A graph $G$ is \emph{$k$-correspondence-colorable} if it is $(L,C)$-colorable for all correspondence assignments with lists of size $k$.
The \emph{correspondence chromatic number} of $G$, denoted $\corres{G}$ after Dvo\v{r}\'{a}k and Postle, is the least $k$ such that $G$ is $k$-correspondence-colorable.

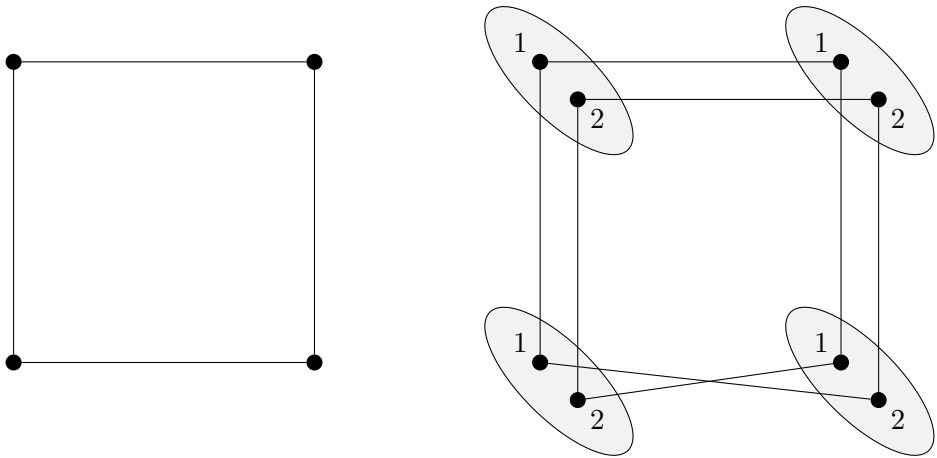
\begin{figure}[ht]
  \centering
    \begin{tikzpicture}[every node/.style={circle,draw=black,inner sep=2pt,fill}]
        \node (A) at (4,4) {};
        \node (B) at (4,0) {};
        \node (C) at (0,0) {};
        \node (D) at (0,4) {};
        \draw (A) -- (B) -- (C) -- (D) -- (A);
    
    \begin{scope}[shift={(7,0)}]
      \draw[rotate around={45:(.25,3.75)},fill=black!5!white] (.25,3.75) ellipse (.5 and 1.3);
      \draw[rotate around={45:(4.25,3.75)},fill=black!5!white] (4.25,3.75) ellipse (.5 and 1.3);
      \draw[rotate around={45:(4.25,-.25)},fill=black!5!white] (4.25,-.25) ellipse (.5 and 1.3);
      \draw[rotate around={45:(.25,-.25)},fill=black!5!white] (.25,-.25) ellipse (.5 and 1.3);
      
      \node (A1) at (0,4) [label=above left:$1$,circle,draw,inner sep=2pt,fill] {};
      \node (A2) at (.5,3.5) [label=below right:$2$,circle,draw,inner sep=2pt,fill] {};
      \node (B1) at (4,4) [label=above left:$1$,circle,draw,inner sep=2pt,fill] {};
      \node (B2) at (4.5,3.5) [label=below right:$2$,circle,draw,inner sep=2pt,fill] {};
      \node (C1) at (4,0) [label=above left:$1$,circle,draw,inner sep=2pt,fill] {};
      \node (C2) at (4.5,-0.5) [label=below right:$2$,circle,draw,inner sep=2pt,fill] {};
      \node (D1) at (0,0) [label=above left:$1$,circle,draw,inner sep=2pt,fill] {};
      \node (D2) at (.5,-.5) [label=below right:$2$,circle,draw,inner sep=2pt,fill] {};
      
      \draw (A1)--(B1)--(C1)--(D2)--(A2);
      \draw (A2)--(B2)--(C2)--(D1)--(A1);
    \end{scope}
  \end{tikzpicture}
  \caption{The cycle $C_4$ on four vertices has Alon--Tarsi number $2$, but is not 2-correspondence-colorable.}
  \label{fig:C_4}
\end{figure}

By the above argument, if a graph $G$ is $k$-correspondence-colorable, then it is also $k$-choosable.
However, the converse does not hold in general.
For example, a cycle on four vertices (see Figure \ref{fig:C_4}), is $2$-choosable, but not $2$-correspondence-colorable. 
This can be seen by assigning a list of two colors, $[1,2]$ to each vertex, then assigning the correspondence which maps each edge except the top one to the matching $1 \to 1$, $2 \to 2$ while the top edge receives the matching $1 \to 2$, $2 \to 1$.

\paragraph{Alon--Tarsi number.}

An \emph{Eulerian subdigraph} $F$ of a digraph $D$ is a spanning subdigraph of $D$ such that for all vertices $v$, $d^+_F(v) = d^-_F(v)$.
An Eulerian subgraph is said to be \emph{even} if it has an even number of edges, and \emph{odd} if it has an odd number of edges.
The \emph{Alon--Tarsi number $\alontarsi{G}$} of $G$ is the minimum $k$ such that there exists an orientation $D$ of $G$ with the outdegree of each vertex less than $k$ and the number of even Eulerian subdigraphs of $D$ differing from the number of odd Eulerian subdigraphs.

Note that the Alon--Tarsi number is finite, since any acyclic orientation of $G$ has one even Eulerian subdigraph (the spanning subdigraph with no edges) and no odd Eulerian subdigraphs.

Alon and Tarsi~\cite{AlonTarsi1992} used the Combinatorial Nullstellensatz to prove that the Alon--Tarsi number is an upper bound on the list chromatic number.
Shauz~\cite{Schauz2010} later gave a proof not using algebraic methods.

To show a lower bound on the Alon--Tarsi number that we use later we need the following lemma.

\begin{lemma}\label{lem:maxdegorientation}
  If a graph $G$ has an orientation $D$ with maximum outdegree $\Delta^+$, then
  \[ \Delta^+ \ge \frac{\mad{G}}{2}. \]
\end{lemma}

\begin{proof}
  Suppose that $D$ is an orientation of $G$ with maximum outdegree $\Delta^+$.
  For any subdigraph $F$ of $D$,
  \[ |E(F)| = \sum_{v\in V(F)} d_F^+(v) \le \sum_{v\in V(F)} d_D^+(v) \le \Delta^+ \cdot |V(F)|, \]
  and hence \[ \Delta^+ \ge \frac{|E(F)|}{|V(F)|}. \]
  As this relation holds for any subdigraph,
  \[ \Delta^+ \ge \frac{\mad{G}}{2}. \qedhere \]
\end{proof}

In fact, the condition given in the previous lemma is also sufficient, which follows from the characterization by Hakimi~\cite{Hakimi1965} of degree sequences of orientations of an undirected graph.

\begin{corollary}\label{cor:AT_bounded_by_mad}
  For a graph $G$,
  \[ \alontarsi{G} > \frac{\mad{G}}{2}. \]
\end{corollary}

\begin{proof}
  Let $D$ be an orientation of $G$ such that the outdegree of each vertex is less than $\alontarsi{G}$ and such that the number of even Eulerian subdigraphs of $D$ differs from the number of odd Eulerian subdigraphs.
  By Lemma~\ref{lem:maxdegorientation}, the maximum outdegree $\Delta^+$ of $D$ is at least $\mad{G}/2$.
  Thus, \[ \alontarsi{G} > \Delta^+ \ge \frac{\mad{G}}{2}. \qedhere \]
\end{proof}

\section{Alon--Tarsi number larger than the correspondence chromatic number}

There exist many graphs for which the Alon--Tarsi number is larger than the correspondence chromatic number, and in fact the ratio between the two parameters can be arbitrarily large.
This follows from a result of Bernshteyn on the asymptotics of the correspondence chromatic number.

\begin{theorem}[Bernshteyn~\cite{bernshteyn2016}]
    \label{thm:berntrifree}
    There exists a positive constant $C$ such that for any triangle-free graph $G$ with maximum degree $\Delta$,
    \[ \corres{G} \le C \frac{\Delta}{\ln \Delta}. \]
\end{theorem}

\begin{corollary}
    There exist infinitely many graphs $G$ where
    \[ \corres{G} < \alontarsi{G}, \]
    and the ratio between the parameters can be arbitrarily large.
\end{corollary}
\begin{proof}
    All $d$-regular graphs have maximum average degree $d$, and so by Corollary~\ref{cor:AT_bounded_by_mad} we have
    \[ \alontarsi{G} > \frac{d}{2}. \]
    All $d$-regular triangle-free graphs have maximum degree $d$, and so by Theorem~\ref{thm:berntrifree} we have
    \[ \corres{G} \le C \frac{d}{\ln d} \]
    for some positive constant $C$.
    For large enough $d$, $C \frac{d}{\ln d} \le \frac{d}{2}$, and the ratio between them becomes arbitrarily large as $d$ grows.
\end{proof}

\section{Correspondence chromatic number larger than the Alon--Tarsi number}

We present in this section a sequence of graphs where the correspondence chromatic number is one larger than the Alon--Tarsi number.
The first graph of the sequence is the cycle on four vertices, and the second graph is a graph on six vertices that we discovered through computation.

\begin{figure}
  \centering
    \begin{subfigure}[b]{0.4\textwidth}
    \centering
    \begin{tikzpicture}
      \node (U) at (0, 2) [label=above:$u$,circle,draw,inner sep=2pt,fill] {};
      \node (V) at (0,-2) [label=below:$v$,circle,draw,inner sep=2pt,fill] {};
      \node[ellipse,draw=black,minimum height=20pt,minimum width=35pt] (K1) at ( 1.5,0) {$K_m$};
      \node[ellipse,draw=black,minimum height=20pt,minimum width=35pt] (K2) at (-1.5,0) {$K_m$};
      \node at (-2.65,.15) {$Q_1$};
      \node at (2.65,.15) {$Q_2$};
      
      \draw[shorten >=5pt] (U) -- (K2.east);
      \draw[shorten >=20pt] (U) -- ($(K2.west) + (-10pt,0)$);
      \draw[shorten >=5pt] (U) -- (K2);
      \draw[shorten <=5pt] (K2.east) -- (V);
      \draw[shorten <=20pt] ($(K2.west) + (-10pt,0)$) -- (V);
      \draw[shorten <=5pt] (K2) -- (V);
      \draw[shorten >=5pt] (V) -- (K1.west);
      \draw[shorten >=20pt] (V) -- ($(K1.east) + (10pt,0)$);
      \draw[shorten >=5pt] (V) -- (K1);
      \draw[shorten <=5pt] (K1.west) -- (U);
      \draw[shorten <=20pt] ($(K1.east) + (10pt,0)$) -- (U);
      \draw[shorten <=5pt] (K1) -- (U);
    \end{tikzpicture}
    \end{subfigure}
    ~ 
    \begin{subfigure}[b]{0.4\textwidth}
    \centering
    \begin{tikzpicture}
      \node (U) at (0, 2) [label=above:$u$,circle,draw,inner sep=2pt,fill] {};
      \node (V) at (0,-2) [label=below:$v$,circle,draw,inner sep=2pt,fill] {};
      \node[ellipse,draw=black,minimum height=20pt,minimum width=35pt] (K1) at ( 1.5,0) {$K_m$};
      \node[ellipse,draw=black,minimum height=20pt,minimum width=35pt] (K2) at (-1.5,0) {$K_m$};
      \node at (-2.65,.15) {$Q_1$};
      \node at (2.65,.15) {$Q_2$};
      
        \draw[-latex, shorten >=5pt] (U) -- (K2.east);
        \draw[-latex, shorten >=20pt] (U) -- ($(K2.west) + (-10pt,0)$);
        \draw[-latex, shorten >=5pt] (U) -- (K2);
        \draw[-latex, shorten <=5pt] (K2.east) -- (V);
        \draw[-latex, shorten <=20pt] ($(K2.west) + (-10pt,0)$) -- (V);
        \draw[-latex, shorten <=5pt] (K2) -- (V);
        \draw[-latex, shorten >=5pt] (V) -- (K1.west);
        \draw[-latex, shorten >=20pt] (V) -- ($(K1.east) + (10pt,0)$);
        \draw[-latex, shorten >=5pt] (V) -- (K1);
        \draw[-latex, shorten <=5pt] (K1.west) -- (U);
        \draw[-latex, shorten <=20pt] ($(K1.east) + (10pt,0)$) -- (U);
        \draw[-latex, shorten <=5pt] (K1) -- (U);
    \end{tikzpicture}
    \end{subfigure}
    \caption{The construction $G_n$ and its orientation $D$.}
    \label{fig:k_large_graph}
\end{figure}
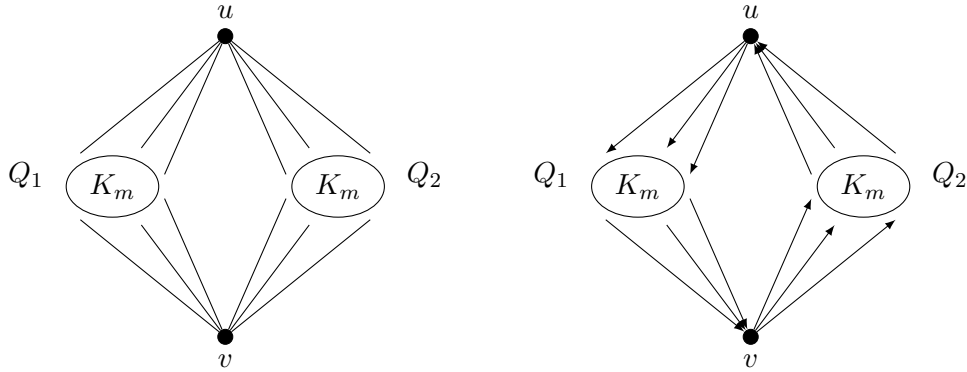

\begin{definition}
  Let $n$ be an even integer at least $4$.
  The $n$-vertex graph $G_n$ is constructed by taking two disjoint graphs isomorphic to $K_m$, where $m = \frac{1}{2}(n-2)$, and adding every possible edge to two independent vertices, $u$ and $v$.
  In other words, $G_n$ is $2K_m$ joined to $2K_1$.
  See Figure~\ref{fig:k_large_graph}.
\end{definition}

We will show that the correspondence chromatic number of $G_n$ is one more than its Alon--Tarsi number.
To prove the Alon--Tarsi number of $G_n$, we consider an appropriate orientation and analyze its Eulerian subdigraphs.
We will need the following definition.

\begin{definition}  
  Given a digraph $D$, an \emph{almost-Eulerian subdigraph $AE_{j, s \to t}$ of order $j$ with source $s$ and sink $t$} is a subdigraph $D'$ of $D$ such that
  \begin{align*}
      d^+_{D'}(s) &= d^-_{D'}(t) = j, \\
      d^-_{D'}(s) &= d^+_{D'}(t) = 0, \\
      d^+_{D'}(v) &= d^-_{D'}(v) \hspace{10pt} \text{for all $v\ne s,t$.}
  \end{align*}
\end{definition}
Note that if $AE_{j, s \to t}$ and $AE_{j, t \to s}$ are edge-disjoint Eulerian subdigraphs, then their union is an Eulerian subdigraph of $D$.

\begin{theorem}
  \label{thm:alontarsi_G}
  Let $n$ be an even integer at least $4$.
  Then \[ \alontarsi{G_n} = \frac{n}{2}. \]
\end{theorem}
\begin{proof}
  Since $G_n$ contains a clique of size $\frac{n}{2}$, $\alontarsi{G_n} \ge \frac{n}{2}$.

  To establish the upper bound, we present an orientation $D$ of $G_n$ such that all vertices have outdegree at most $m = \frac{n}{2} - 1$, and such that the number of even Eulerian subdigraphs differs from the number of odd Eulerian subdigraphs.
  
  To form $D$, we orient the edges from $u$ towards the clique $Q_1$, and we orient the edges from $Q_2$ towards $u$.
  We orient the edges incident to $v$ oppositely:
  we orient the edges from $Q_1$ towards $v$, and we orient the edges from $v$ towards the clique $Q_2$.
  We orient the edges within the cliques $Q_1$ and $Q_2$ acyclically, so that they are isomorphic to the transitive tournament on $m$ vertices.
  See Figure~\ref{fig:k_large_graph} for an illustration of $D$.
  Observe that the maximum outdegree of the orientation $D$ is $m$.
  
  Let $F$ be an Eulerian subdigraph of $D$.
  We can decompose $F$ into two pieces: the edges $X$ outgoing from $u$, incident to vertices of $Q_1$, or incoming to $v$; and the edges $Y$ outgoing from $v$, incident to vertices of $Q_2$, or incoming to $u$.
  Note that $X$ (thought of as a spanning subdigraph of $D$ with the specified edges) is an almost-Eulerian subdigraph $AE_{j, u \to v}$, where $j$ is the number of outgoing edges from $u$ (and also the number of incoming edges to $v$), and similarly for $Y$.
  We can reverse this process to count the Eulerian subdigraphs: the union of any almost-Eulerian subdigraph $AE_{j, u \to v}$ and any almost-Eulerian subdigraph $AE_{j, u \to v}$ is an Eulerian subdigraph of $D$.
  
  Let $E_j$ be the set of almost Eulerian subdigraphs $AE_{j, u \to v}$ of $D$ with an even number of edges, and let $O_j$ be the set of almost Eulerian subdigraphs $AE_{j, u \to v}$ with an odd number of edges.
  Similarly, let $E'_j$ and $O'_j$ be the sets of even and odd, respectively, almost Eulerian subdigraphs $AE_{j, v \to u}$.
  Note that every subdigraph of $E_j \cup O_j$ is edge disjoint from every subdigraph of $E'_j \cup O'_j$.
  Hence the union of any $X \in E_j \cup O_j$ and $X' \in E'_j \cup O'_j$ is an Eulerian subdigraph of the orientation $D$ of $G_n$.
  Since removing $u$ or $v$ from $D$ leaves an acyclic digraph, all Eulerian subdigraphs are constructed in this way.
  
  Note that because of the symmetry of the orientation, $|E_j| = |E'_j|$ and $|O_j| = |O'_j|$. Let $e_j = |E_j| = |E'_j|$ and let $o_j = |O_j| = |O'_j|$.
  The total number of even Eulerian subdigraphs is
  \[ \sum_{j = 0}^m (e_j^2 + o_j^2). \]
  The total number of odd Eulerian subdigraphs is
  \[ \sum_{j = 0}^m 2e_jo_j. \]
  Therefore, their difference is
  \begin{align*}
      \sum_{j=0}^m (e_j^2 + o_j^2) - \sum_{j=0}^m 2e_jo_j &= \sum_{j=0}^m (e_j^2 + o_j^2 - 2e_jo_j) \\
      &= \sum_{j=0}^m (e_j - o_j)^2.
  \end{align*}
  Since this expression is a sum of squares of real numbers, it will be positive unless all the squares are equal to $0$.
  However, $E_0$ contains the empty subdigraph with no edges, while $O_0$ is empty, and so $e_0 - o_0 \ge 1$.
  Hence, this sum of squares is positive. 
  Therefore, the number of even Eulerian subdigraphs differs from the number of odd Eulerian subdigraphs in this orientation, and so $\alontarsi{G_n} \le m + 1 = \frac{n}{2}$.
\end{proof}

\begin{theorem}
  \label{thm:corres_G}
  Let $n$ be an even integer at least $4$.
  Then \[ \corres{G_n} = \frac{n}{2} + 1. \]
\end{theorem}
\begin{proof}
    Since every subgraph of $G_n$ has a vertex of degree at most $m+1 = \frac{n}{2}$, we see that by degeneracy $\corres{G_n} \le n/2 + 1$.
    
    We establish a lower bound by finding a correspondence of $G_n$ with $k = \frac{n}{2}$ colors which cannot be colored.
    Note that the cliques $Q_1$ and $Q_2$ in $G_n$ are of size $k-1$.
    
    Let the correspondences on every edge be the identity correspondence, mapping color $i$ to color $i$, except for on the edges from $u$ to the clique $Q_2$.
    On those edges, the correspondence will map color $i$ on $u$ to color $i+1$ on the clique, and cyclically map color $k-1$ in $u$ to color $0$ in the clique.
    
    Suppose that $u$ is colored $0$ (other cases are similar through cyclic permutation of the colors).
    Then the vertices of $Q_1$ must use all of the colors $1, \ldots, k-1$ in some order. 
    Similarly, the vertices of $H_2$ are forced to use all of the colors $0$ and $2, \ldots, k-1$ in some order.
    Therefore, every color is in the neighborhood of $v$, and hence $v$ cannot be colored. Thus, $\corres{G_n} > n/2$.
\end{proof}

The previous theorems raise the question of the relationship between the correspondence chromatic number and the Alon--Tarsi number when the correspondence chromatic number is larger.
Unlike the situation when the Alon--Tarsi number is larger, the following theorem shows that, in fact, the ratio between the two parameters is bounded.

\begin{theorem}
    \label{thm:corres_alontarsi_factor2}
    Let $G$ be a graph.  Then
    \[ \corres{G} \le 2 \alontarsi{G}. \]
\end{theorem}
\begin{proof}
  Note that the degeneracy of $G$ is less than the maximum average degree $\mad{G}$.
  Therefore,
    \[ \corres{G} \le \mad{G} + 1. \]
  By Corollary~\ref{cor:AT_bounded_by_mad}, 
    \[ \mad{G} < 2\alontarsi{G}. \]
  Combining these, we obtain
    \[ \corres{G} < 2\alontarsi{G} + 1, \]
  from which the result follows.
\end{proof}

\section{Conclusion}

Note that the graphs $G_n$ have the property that their Alon--Tarsi number is only one less than their correspondence chromatic number.
Theorem~\ref{thm:corres_alontarsi_factor2} allows for an arbitrarily large difference between these two quantities, and yet we have not been able to find a graph with a difference of two between these quantities. This leads to the following question.

\begin{question}
  Are there graphs $G$ such that $\alontarsi{G} + 2 \le \corres{G}$?
\end{question}

\bibliographystyle{abbrv}
\bibliography{citations}

\end{document}